# Stochastic integrals and asymptotic analysis of canonical von Mises statistics based on dependent observations

**Igor S. Borisov[1,*] and Alexander A. Bystrov[1,†]**

*Sobolev Institute of Mathematics*

**Abstract:** In the first part of the paper we study stochastic integrals of a nonrandom function with respect to a nonorthogonal Hilbert noise defined on a semiring of subsets of an arbitrary nonempty set.

In the second part we apply this construction to study limit behavior of canonical (i.e., degenerate) Von Mises statistics based on weakly dependent stationary observations.

## 1. Stochastic integrals of non-random kernels for non-orthogonal noises.

### *1.1. Introduction and statement of the main result.*

Let $\{\Omega, \Theta, \mathbf{P}\}$ be a probability space, $\mathfrak{X}$ be an arbitrary nonempty set, and $\mathfrak{M}$ be a semiring with identity of its subsets (i.e., $\mathfrak{X} \in \mathfrak{M}$ and, for all $A, B \in \mathfrak{M}$, we have $A \cap B \in \mathfrak{M}$ and $A \setminus B = \sum_{i \leq n} C_i$, where $C_i \in \mathfrak{M}$). We call a random process $\{\mu(A),\ A \in \mathfrak{M}\}$ *an elementary stochastic measure or a noise* if $\mu(A) \in L_2(\Omega, \mathbf{P})$ for all $A \in \mathfrak{M}$ (i.e., $\mu(\cdot)$ is a *Hilbert process*) and (N1) $\quad \mu(A_1 \cup A_2) = \mu(A_1) + \mu(A_2)$ a.s. if only $A_1 \cap A_2 = \emptyset$ and $A_1 \cup A_2 \in \mathfrak{M}$.

A noise $\mu$ is called *orthogonal* if
(N2) $\quad \mathbf{E}\mu(A_1)\mu(A_2) = m_0(A_1 \cap A_2),$
where $m_0$ a finite measure (*the structure function*) on $\sigma(\mathfrak{M})$ [14].

**Typical Examples.**

(i) Consider the following semiring of subsets of a closed interval $[0, T]$:

$$\mathfrak{M} = \{(t, t+\delta];\ 0 < t < t+\delta \leq T\} \bigcup \{[0, \delta];\ 0 < \delta \leq T\}.$$

A random process $\xi(t)$ defined on $[0, T]$ generates the noise

$$\mu((t, t+\delta]) := \xi(t+\delta) - \xi(t),$$

where, in the case $t = 0$, this formula defines the measure of the closed interval $[0, \delta]$. If $\xi(t)$ is a process with independent increments then $\mu$ is an orthogonal noise.

[1]630090, Novosibirsk, Russia, acad. Koptyug avenue, 4, Sobolev Institute of Mathematics, e-mail: sibam@math.nsc.ru; bystrov@ngs.ru

[*]Supported in part by the Russian Foundation for Basic Research Grants 05-01-00810 and 06-01-00738 and INTAS Grant 03-51-5018.

[†]Supported in part by the Russian Foundation for Basic Research Grant 06-01-00738.

(ii) To construct multiple stochastic integrals the semiring $\mathfrak{M}^k := \mathfrak{M} \times \cdots \times \mathfrak{M}$ is considered, where $\mathfrak{M}$ is defined above, and the following *multiple noise* is defined by increments of a random process $\xi(t)$:

$$\mu((\overline{t}, \overline{t} + \overline{\delta}]) = \prod_{i \leq k}(\xi(t_i + \delta_i) - \xi(t_i)),$$

where $(\overline{t}, \overline{t} + \overline{\delta}] = (t_1, t_1 + \delta_1] \times \cdots \times (t_k, t_k + \delta_k]$.

It worth noting that, in general, in the second example the noise $\mu$ does not satisfy condition (N2) even if the process $\xi(t)$ has independent increments.

We note some significant results in the area under consideration:

**I.** *Univariate stochastic integrals based on orthogonal noises*:
   N. Wiener, 1923.
   A. N. Kolmogorov, H. Cramér, 1940.
   I. I. Gikhman and A. V. Skorokhod, 1977.

**II.** *Multiple stochastic integral with a multiple noise generated by a process with independent increments*:
   N. Wiener, 1938, 1958.
   K. Itô, 1951.
   P. Major, 1981.

**III.** *Univariate stochastic integral with a noise generated by increments of a Hilbert process on the real line (nonorthogonal noise)*:
   M. Loève, 1960.
   S. Cambanis and S. Huang, 1978.
   V. Pipiras and M. S. Taqqu, 2000.

**IV.** *Multiple stochastic integral with a multiple noise generated by increments of a Gaussian process on the real line (nonorthogonal noise)*:
   S. Cambanis and S. Huang, 1978.
   A. Dasgupta and G. Kallianpur, 1999.

**General Case.** We begin to study stochastic integrals with nonorthogonal noises defined on semirings of subsets of an arbitrary measurable space. We follow the generality considered in [14], where general stochastic integrals with orthogonal noises were studied. Complete proofs of some statements in the first Section of the paper are published in [4].

Introduce the function $m(A \times B) := E\mu(A)\mu(B)$ indexed by elements of $\mathfrak{M}^2$.

**Main Assumption.** The function $m$ is a finite $\sigma$-additive signed measure (covariance measure) on $\mathfrak{M}^2$.

**Example.** Let $\Phi(t, s) := \mathbf{E}\xi(t)\xi(s)$ be the covariance function of a centered Hilbert random process $\xi(t)$ defined on a closed interval $[0, T]$. We say that the function $\Phi(t, s)$ possesses a bounded variation if, for a constant $C$,

$$\sup_{\{t_i, s_i\}} \sum_{i,j} |\Delta\Phi(t_i, s_j)| \leq C,$$

where $\Delta\Phi(t_i, s_i) := \Phi(t_{i+1}, s_{j+1}) + \Phi(t_i, s_j) - \Phi(t_{i+1}, s_j) - \Phi(t_i, s_{j+1})$ (*the double difference*); $0 = t_0 < t_1 < \cdots < t_n = T$, $0 = s_0 < s_1 < \cdots < s_l = T$ are arbitrary finite partitions of the interval $[0, T]$, and the supremum is taken over all such partitions (see also [9]).

**Proposition 1.** *If $\Phi(t,s)$ has a bounded variation then the Main Assumption for the corresponding covariance measure is valid*

It is well known that any measure of such a kind can be uniquely extended onto $\sigma(\mathfrak{M}^2)$. Moreover, $m = m^+ - m^-$ (Hahn–Jordan decomposition), where $m^+$ and $m^-$ are nonnegative finite measures. Put $|m| := m^+ + m^-$ (*the total variation measure*).

Introduce the space of $\sigma(\mathfrak{M})$-measurable functions:

$$S := \{f : \int_{\mathfrak{X}^2} f(t)f(s)m(dt \times ds) < \infty\}.$$

For any $\sigma(\mathfrak{M})$-measurable functions $f, g \in S$ consider the bilinear symmetric functional

$$d(f,g) := \int_{\mathfrak{X}^2} f(t)g(s)m(dt \times ds).$$

It is clear that $d(f,f) \geq 0$. But, in general, the equation $d(f,f) = 0$ has not only zero solution. Denote $\|f\| := d(f,f)^{1/2}$ (seminorm). If $S$ is the factor space w.r.t. to the condition $d(f,f) = 0$ then $S$ is an Euclidean space. But, in general, $S$ may be incomplete (i.e., is not Hilbert) ([22]).

Notice that the space $S$ can also defined in such a way:

$$S = \{f : \int_{\mathfrak{X}^2} |f(t)f(s)||m|(dt \times ds) < \infty\}.$$

But the functional $\|f\|_* := \left(\int_{\mathfrak{X}^2} |f(t)f(s)||m|(dt \times ds)\right)^{1/2}$ may not satisfy the triangle inequality ($\|f\|_*$ is a seminorm iff $|m|$ is *nonnegatively defined*).

For an orthogonal noise $\mu$ with a structure function $m_0$ the following obvious equality chain is valid: $\|f\| = \|f\|_* = \|f\|_{L_2(\mathfrak{X}, m_0)}$.

**Theorem 1.** *Let $f \in (S, \|\cdot\|)$. Then there exists a sequence of step functions*

(1) $$f_n(x) := \sum_{k \leq n} c_k I(x \in A_k),$$

*where $c_k \in \mathbb{R}, A_k \in \mathfrak{M}$, converging in $(S, \|\cdot\|)$ to $f$ as $n \to \infty$. Moreover, the sequence*

(2) $$\eta(f_n) := \sum_{k \leq n} c_k \mu(A_k)$$

*mean-square converges to a limit random variable $\eta(f)$ which does not depend on the sequence of step functions $\{f_n\}$.*

*Proof.* (For detail see [4]). Let $\{f_n\}$ be a sequence of step functions of the form (1) converging to $f$ in the seminorm $\|\cdot\|$. One can prove that the sequences of such a kind exist. Then

$$\|\eta(f_n) - \eta(f_k)\|^2_{L_2(\Omega, \mathbf{P})} = \sum_{i,j}(c_i^{(n)} - c_i^{(k)})(c_j^{(n)} - c_j^{(k)})m(A_i \times A_j)$$

$$\equiv \|f_n - f_k\|^2 \to 0$$

as $n, k \to \infty$ due to the triangle inequality (i. e., $\{f_n\}$ is a Cauchy sequence in $S$). Without loss of generality, we may assume here that the step functions $f_k$ and $f_n$ are defined on a common partition $\{A_i\}$. Hence $\{\eta(f_n)\}$ is a Cauchy sequence in the Hilbert space $L_2(\Omega, \mathbf{P})$. Thus, in this Hilbert space, there exists a limit random variable $\eta(f)$ for the sequence $\{\eta(f_n)\}$. □

**Remark 1.** Since, in general, the space $S$ is incomplete then, in this case, one cannot construct an *isometry* (one-to-one mapping preserving distances) between the $L_2$-closed linear span of all integral sums and $S$ ([22]). Existence of such isometry is a key argument of the classical construction of stochastic integrals with orthogonal noises.

**Remark 2.** The generality in Theorem 1 allows us to define stochastic integrals both for the univariate and multivariate cases studied by the predecessors mentioned above. In the case of Gaussian noises generated by arbitrary centered Gaussian processes on the real line, our construction differs from that in [9], where multiple stochastic integrals are defined by the corresponding tensor power of the reproducing Hilbert space corresponding to the above-mentioned initial Gaussian process. However, one can prove that, to define these multiple integrals in the Gaussian case, the descriptions of the corresponding kernel spaces in these two constructions coincide.

**Remark 3.** If we consider the introduced-above multiple stochastic integral with the product–noise generated by a White noise with a structure function $m_0$, and, moreover, the kernel vanishes on all diagonal subspaces then our construction coincides with the classical Wiener – Itô multiple construction. Notice that, in this construction, for the kernels with zero values on all diagonal subspaces, there exists the isometry mentioned in Remark 1. In this case, the space $S$ coincides with the Hilbert space $L_2(\mathfrak{X}^k, m_0^k)$, where $k$ is the dimension of the multiple integral.

### 1.2. Infinitesimal analysis of covariance measures.

We now describe some function kernel spaces to define the stochastic integrals in Theorem 1. We start with the univariate construction.

#### 1.2.1. Univariate stochastic integral.

Consider a centered random process $\xi(t)$ with a covariance function $\Phi(t, s)$. In all the examples of Section I we put $\mathfrak{X} = [0, T]$. This process generates the elementary stochastic measure $\mu(dt) := d\xi(t)$ introduced above.

**Regular covariance functions.** In the above-mentioned definition of the double difference of the covariance function $\Phi$ we set $t_{j+1} := t_j + \delta$ and $s_{j+1} := s_j + \delta$. Assume that, for all $\delta > 0$ and $t_j, s_j, t_j + \delta, s_j + \delta \in [0, T]$,

$$\Delta \Phi(t_i, s_i) = \int_{t_i}^{t_i+\delta} \int_{s_i}^{s_i+\delta} q(t,s) \lambda(dt) \lambda(ds),$$

where $\lambda$ is an arbitrary $\sigma$-finite measure. If

$$\int_{\mathfrak{X}^2} |f(t)f(s)q(t,s)| \lambda(dt)\lambda(ds) < \infty$$

then $f \in S$ (i.e., $\int_{\mathfrak{X}} f(t) d\xi(t)$ is well defined.)

For example, the regular FBM has the covariance function

$$\Phi(t,s) = \frac{1}{2}(t^{2h} + s^{2h} - |t-s|^{2h}),$$

where $h \in (1/2, 1]$. In this case, $q(t,s) := h(2h-1)|t-s|^{2h-2}$, $t \neq s$, and $\lambda(dt) := dt$ is the Lebesgue measure. Moreover, in this case one can prove the embedding $L_{1/h}(\mathfrak{X}, dt) \subseteq S$. ([22].)

**Irregular covariance functions.** Consider the class of *factorizing* covariance functions

$$\Phi(t,s) = G(\min(t,s))H(\max(t,s)).$$

It is known [3] that $\Phi(t,s)$ of such a kind is the covariance function of a nondegenerate on $(0, T)$ random process iff the fraction $\frac{G(t)}{H(t)}$ is a nondecreasing positive function. In particular, any Gaussian Markov process with non-zero covariance function admits such factorization: For example, a standard Wiener process has the components $G(t) = t$ and $H(t) \equiv 1$; a Brownian bridge on $[0,1]$ has the components $G(t) = t$ and $H(t) = 1 - t$; and, finally, an arbitrary stationary Gaussian process on the positive half-line has the components $G(t) = \exp(\alpha t)$ and $H(t) = \exp(-\alpha t)$, where $\alpha > 0$.

Notice that, in this case, the function $\Phi(t,s)$ is nondifferentiable on the diagonal if the components are nondegenerate functions.

Let, in addition, $G(t) \uparrow$, $H(t) \downarrow$ be monotone, positive on $(0, T)$, and absolutely continuous w.r.t. the Lebesgue measure on $[0, T]$. We prove that $supp\, m^- = [0, T]^2 \setminus D$, where $D := \{(t,s) : t = s\}$, is the main diagonal of the square $[0, T]^2$. Let $s < s + \delta \leq t < t + \Delta$. Then

$$\begin{aligned} m(\,(s, s+\delta] &\times (t, t+\Delta]\,) \\ &= G(s)H(t) + G(s+\delta)H(t+\Delta) - G(s)H(t+\Delta) - G(s+\delta)H(t) \\ &= (G(s+\delta) - G(s))(H(t+\Delta) - H(t)) < 0. \end{aligned}$$

In other words, $m^-$ is absolutely continuous w.r.t. the bivariate Lebesgue measure $\lambda_2$ and the corresponding Radon – Nikodym derivative is defined by the formula

$$\frac{dm^-}{d\lambda_2}(t,s) = G'(t)|H'(s)|.$$

We now calculate $m^+$-measure of an infinitesimal diagonal square:

$$\begin{aligned} m^+((s, s+h] &\times (s, s+h]) \\ &= H(s+h)\bigl(G(s+h) - G(s)\bigr) - G(s)\bigl(H(s+h) - H(s)\bigr) \\ &= \int_s^{s+h} \bigl(H(z)G'(z) - G(z)H'(z)\bigr)dz + o(h). \end{aligned}$$

Hence $m^+$ is absolutely continuous w.r.t. the induced Lebesgue measure on the diagonal.

Therefore, in the case under consideration, the measures $m^-$ and $m^+$ are *singular*. Finally, to verify the condition $f \in S$ we need to verify existence of the following two integrals:

$$\int_{\mathfrak{X}} f(t)^2 (H(t)+1) dG(t), \quad \int_{\mathfrak{X}} f(t)^2 (G(t)+1) dH(t).$$

**Covariance functions of mixed type.** Let $\{X_n;\ n \geq 1\}$ be a stationary sequence of r.v.'s satisfying $\varphi$-mixing condition. Consider a centered Hilbert process (not necessarily Gaussian!) $Y(t)$ with the covariance function which is well defined under some restrictions on the $\varphi$-mixing coefficient (see the Gaussian case in [1]):

$$\mathbf{E}Y(s)Y(t) = F(\min(s,t)) - F(t)F(s) + \sum_{j \geq 1} \left(F_j(s,t) + F_j(t,s) - 2F(s)F(t)\right),$$

where $F(t)$ is the distribution function of $X_1$ which is assumed absolutely continuous with a density $p(t)$, and $F_j(t,s)$ is the joint distribution functions of the pairs $(X_1, X_{j+1})$. Let all the functions $F_j(t,s)$, $j = 1, 2, \ldots$, have bounded densities $p_j(t,s)$. For all $t, s \in R$, we assume that the series

$$b(t,s) := \sum_{j \geq 1} [p_j(t,s) + p_j(s,t) - 2p(t)p(s)]$$

absolutely converges and the corresponding series $\sum |\cdot|$ is integrable on $R^2$.

We note that, under the above-mentioned restrictions, we deal with a covariance function represented as a sum of covariance functions from items A.1 and A.2. Hence we may use the infinitesimal analysis of the corresponding covariance measures from these items. Indeed,

$$m((t, t + \Delta] \times (s, s + \Delta])$$
$$= \mathbf{P}(X_1 \in (t, t+\Delta] \cap (s, s+\Delta]) + \int_t^{t+\Delta} \int_s^{s+\Delta} (b(u,v) - p(u)p(v))\, du dv.$$

So, under the conditions $\int f^2(t)p(t)dt < \infty$ and $\int \int |f(t)f(s)b(t,s)| dtds < \infty$, we can correctly define the stochastic integral $\int f(t) dY(t)$.

*1.2.2. Multiple stochastic integral.*

We study multiple stochastic integrals (MSI) based on a product-noise defined in Example (**ii**) by increments of a *Gaussian* process. In this case, to calculate the covariance measure we use the following well-known convenient representation:

$$m((t_1, t_1 + \delta] \times \cdots \times (t_{2k}, t_{2k} + \delta]) = \sum \prod \Delta \Phi(t_i, t_j),$$

where the sum is taken over all partitions on pairs of the set $\{1, 2, \ldots, 2k\}$, and the product is taken over all pairs in a fixed such partition. Notice that, in the sequel, to define multiple stochastic integrals we use only this property of Gaussian processes. However, we may study the multiple integrals for the non-Gaussian case: For example, if the integrating process $\xi(t)$ can be represented as a polynomial transform of a Gaussian process. In this case, to define the multiple integral, we can obtain some restrictions on the kernel $f$ close to the conditions below.

**Remark 4.** The Main Assumption for the covariance measure $m(A_1 \times \cdots \times A_{2k})$ introduced above follows from the Proposition 1 if only $\Phi(t,s)$ has a bounded variation. This property of the covariance function is fulfilled in items B.1 – B.3 below.

**Regular covariance functions.** *Conditions to define MSI*:

$$\int_{[0,T]^{2k}} |f(t_1,\ldots,t_k)f(t_{k+1},\ldots,t_{2k})| \sum \prod |q(t_i,t_j)| dt_1 \cdots dt_{2k} < \infty,$$

where the sum and the product are introduced above, $q(t,s)$ is the density (the Radon–Nikodym derivative) of $\Delta\Phi$. As a consequence, we obtain the main result in [7], for the regular FBM from item A.1. In this case we should set in this condition $q(t,s) := h(2h-1)|t-s|^{2h-2}$.

**Factorizing covariance functions.** Let the factorizing components $H$ and $G$ be smooth functions.
1) If $t_i = t_j$ then, as $\delta \to 0$, we have the following asymptotic representation of the double difference of $\Phi(\cdot)$ on the infinitesimal cube $(t_i, t_i+\delta] \times (t_j, t_j+\delta]$:

$$\Delta\Phi(t_i,t_j) = \delta(H(t_i)G'(t_i) - G(t_i)H'(t_i)) + O(\delta^2),$$

2) If $t_i \neq t_j$ then, as $\delta \to 0$,

$$\Delta\Phi(t_i,t_j) = \delta^2 G'(\min(t_i,t_j))H'(\max(t_i,t_j)) + o(\delta^2).$$

Denote
$$g_1(t) := H(t)G'(t) - G(t)H'(t),$$
$$g_2(t,s) := G'(\min(t,s))H'(\max(t,s)).$$

A set $D_{(r_1,\ldots,r_l)} \subset [0,T]^{2k}$ is called a *diagonal subspace determined by variables of multiplicity* $r_1,\ldots,r_l$ ($r_i \geq 2$, $\sum r_i < 2k$) if it defines by the following $l$ chains of equalities:
$$x_{i_{j,1}} = \cdots = x_{i_{j,r_j}} \quad j=1,\ldots,l,$$
where $i_{j,m} \neq i_{n,d}$ for any $(j,m) \neq (n,d)$.

**Proposition 2** (see Borisov and Bystrov, 2006a). *In the case under consideration any covariance measure $m$ has zero mass on any diagonal subspace $D_{(r_1,\ldots,r_l)}$ having at least one multiplicity $r_i > 2$.*

Given a kernel $f(t_1,\ldots,t_k)$ we set $\varphi_f(t_1 \cdots t_{2k}) := f(t_1,\ldots,t_k)f(t_{k+1},\ldots,t_{2k})$.
*Conditions to define MSI:* First, we need to verify the condition

$$(3) \quad \begin{aligned} &\int |\varphi_f(s_1,s_1,\ldots s_n,s_n,t_1,\ldots t_{2(k-n)})| \\ &\times \prod_{i=1}^n g_1(s_i) \sum \prod |g_2(t_i,t_j)| ds_1 \cdots ds_n dt_1 \cdots dt_{2(k-n)} < \infty \end{aligned}$$

for all $n = 0, 1, \ldots, k$, and, second, to verify finiteness of all the analogous integrals for all permutations of $2k$ arguments of the kernel $\varphi_f$. Here, by definition, $\prod_{i=1}^{0} = 1$.

If the kernel $f(\cdot)$ is symmetric and vanishes on all diagonal subspaces, and the functions $g_1$ and $g_2$ are bounded then condition (3) is reduced to the restriction $f \in L_2(\mathfrak{X}^k, dt_1 \cdots dt_k)$. In particular, if the multiple noise is defined by increments of a standard Wiener process then $g_1 \equiv 1$ and $g_2 \equiv 0$ (cf. [17, 19]).

**Covariance functions of mixed type.** We now define the multiple stochastic integral for a Gaussian process $Y(t)$ with the covariance introduced in A.3. Let $p(t)$ and $b(t,s)$ defined in A.3 be continuous functions. Then, as $\delta \to 0$, we have for $t_i = t_j$ (see A.3):
$$\Delta\Phi(t_i, t_j) = \delta p(t_i) + o(\delta).$$

If $t_i \neq t_j$ then
$$\Delta\Phi(t_i, t_j) = \delta^2 \left(b(t_i, t_j) - p(t_i)p(t_j)\right) + o(\delta^2).$$

So, we actually repeat the arguments from item B.2 and to define the multiple stochastic integral for $Y(t)$ we need to verify condition (3) for $g_1(t) := p(t)$ and $g_2(t,s) := b(t,s)$.

## 2. Asymptotics of canonical von mises statistics.

In this Section we consider some applications of the MSI construction from Section I. We study limit behavior of multivariate Von Mises functionals of empirical processes based on samples from a stationary sequence of observations.

Let $\{X_n;\ n \geq 1\}$ be a stationary sequence of $[0,1]$-uniformly distributed r.v.'s satisfying the $\psi$-mixing condition: $\psi(m) \to 0$ if $m \to \infty$, where

$$(4) \qquad \psi(m) := \sup \left| \frac{\mathbf{P}(AB)}{\mathbf{P}(A)\mathbf{P}(B)} - 1 \right|, \quad m = 1, 2, \ldots,$$

and the supremum is taken over all events $A$ and $B$ (having non-zero probabilities) from the respective $\sigma$-fields $\mathcal{F}_1^k$ and $\mathcal{F}_{k+m}^\infty$, where $\mathcal{F}_l^k$, $l \leq k$, is the $\sigma$-field generated by the random variables $X_l, \ldots, X_k$, as well as over all natural $k$. This mixing condition was introduced in Blum, Hanson and Koopmans, 1963.

Introduce the normalized $d$-variate Von Mises statistics (or $V$-*statistics*)

$$(5) \qquad V_n := n^{-d/2} \sum_{1 \leq i_1, \ldots, i_d \leq n} f(X_{i_1}, \ldots, X_{i_d}), \quad n = 1, 2, \ldots,$$

where the kernel $f(\cdot)$ satisfies the degeneracy condition

$$\mathbf{E} f(t_1, \ldots, t_{k-1}, X_k, t_{k+1}, \ldots, t_d) = 0$$

for all $t_1, \ldots, t_d \in [0,1]$ and $k = 1, \ldots, d$. Such *canonical* statistics were introduced in [25], and [15], where, moreover, the so-called $U$-*statistics* were studied:

$$U_n := (C_n^d)^{-1/2} \sum_{1 \leq i_1 < \cdots < i_d \leq n} f_0(X_{i_1}, \ldots, X_{i_d}),$$

where, as a rule, the kernel $f_0$ is symmetric w.r.t. all permutations of the arguments.

Notice that in the definitions above we may consider the observations $\{X_i\}$ taking values in an arbitrary measurable space and having an arbitrary distribution. If the sample distribution has no atoms then an $U$-statistics with a symmetric kernel can be easily reduced to a $V$-statistics with a new kernel having zero values on all diagonal subspaces. For the IID samples limit behavior of these statistics is well known (for reference see Korolyuk and Borovskikh, 1994). In the case of one-dimensional observations, using the corresponding quantile transforms, we may reduce these statistics to those based on samples from the $[0,1]$-uniform distribution.

Introduce the normalized empirical process $S_n(t) := \sqrt{n}(F_n^*(t) - t)$, where $F_n^*(t)$ is the standard empirical distribution function based on the above-introduced sample $X_1, \ldots, X_n$ with the $[0,1]$-uniform marginal distribution. We now recall the well-known key representation of the canonical $V$-statistics (5) by the Lebesgue integral:

$$
V_n = \int_{[0,1]^d} f(x_1, \ldots, x_d) \, dS_n(x_1) \cdots dS_n(x_d). \tag{6}
$$

We recall that a diagonal subspace in $[0,1]^d$ is determined as

$$
D_{i^*,q^*} := \{(t_1, \ldots, t_d) \in [0,1]^d :
$$
$$
t_{i_1} = t_{i_2} \cdots = t_{i_{q_1}}, \ldots, t_{i_{q_{r-1}+1}} = t_{i_{q_{r-1}+2}} \cdots = t_{i_{q_r}} \},
$$

where $i^* := (i_1, i_2, \ldots, i_{q_1}, \ldots, i_{q_{r-1}+1}, i_{q_{r-1}+2}, \ldots, i_{q_r})$ is a vector with pairwise different integer coordinates from the set $\{1, \ldots, d\}$ and $q^* := (q_1, \ldots, q_r)$, where $q_1 \geq 2$, $q_r \leq d$ and $q_{i+1} - q_i \geq 2$ if $1 \leq i < r$. In the sequel, to indicate diagonal subspaces, we will use natural parameter $i$ instead of vector-valued $(i^*, q^*)$ renumbering all the diagonal subspaces.

The subspace of $[0,1]^d$ defined by pairwise different coordinates $t_i$ is called the *main subspace* and denoted by $D_0$. Obviously,

$$
D_0 = [0,1]^d \setminus \bigcup_{i \geq 1} D_i.
$$

Introduce the following function space

$$
S_0 := \{f : \sum_{i \geq 0} \int_{D_i} f^2(t_1, \ldots, t_d) dt_{i_{q_1}} \cdots dt_{i_{q_r}} < \infty\},
$$

where the subscripts $i_{q_1}, \ldots, i_{q_r}$ determine the corresponding diagonal (or the main) subspace $D_i$ of dimension $r$ and the kernel $f^2(\cdot)$ in the corresponding multiple integral $\int_{D_i}$ of this notation has only $r$ independent arguments. For example, if $D_i = \{(t_1, t_2, t_3, t_4) \in [0,1]^4 : t_1 = t_3, \, t_2 = t_4\}$ then

$$
\int_{D_i} f^2(t_1, t_2, t_3, t_4) \, dt_3 dt_4 = \int_{[0,1]^2} f^2(t_3, t_4, t_3, t_4) \, dt_3 dt_4.
$$

Notice that, under the conditions of Theorem 2 below, the function space $S_0$ is a subspace of $S$ introduced in Theorem 1 for $d$-fold product-noise generated by increments of the stochastic processes $Y(t)$ from items A.3 and B.3 since the series

$$
b(t,s) = \sum_{k \geq 1} (p_k(t,s) + p_k(s,t) - 2)
$$

from item A.3 satisfies the restrictions we need and, moreover, it is uniformly bounded on $[0,1]^2$ due to $\psi$-mixing condition (4) and the restrictions on the mixing coefficient in Theorem 2 below. First of all, we observe that (4) implies existence of joint densities $p_k(t,s)$ since we may consider in (4) 'infinitesimal events" $A := \{X_1 \in (t, t+dt)\}$ and $B := \{X_{k+1} \in (s, s+ds)\}$ and take the relations $\mathbf{P}(A) = dt$ and $\mathbf{P}(B) = ds$ into account. Moreover, we obtain the upper bound $|b(t,s)| \leq \sum_{k \geq 1} \psi(k)$. This fact allow us to estimate the covariance measure of the corresponding product-noise generated by increments of $Y(t)$.

Define on the space $S_0$ the combined $L_2$-norm

$$||f||_0^2 := \sum_{i \geq 0} \int_{D_i} f^2(t_1, \ldots, t_d) dt_{i_{q_1}} \cdots dt_{i_{q_r}}, \tag{7}$$

where in the case $i = 0$ (i. e., for the main subspace $D_0$,) we put $i_{q_k} = k$ and $r = d$ in the corresponding summand on the right-hand side of (7). Notice that this norm is stronger than $||\cdot||$ since in the case under consideration we have $g_1(t) := p(t) \equiv 1$ on $[0,1]$ and $g_2(t,s) := b(t,s)$ is a bounded function as well (see B.3 and the upper bounds in (20) and (21) below). In other words, in the case under considerations, the linear normed space $(S_0, ||\cdot||_0)$ is embedded to the corresponding space $(S, ||\cdot||)$ in Theorem 1.

**Theorem 2.** *Under the above-mentioned restrictions let $f \in S_0$ and*

$$\Psi(d) := \sum_{k \geq 1} \psi(k) k^{2d-2} < \infty. \tag{8}$$

*Then, as $n \to \infty$,*

$$V_n \xrightarrow{d} \int_{[0,1]^d} f(t_1, \ldots, t_d) dY(t_1) \cdots dY(t_d),$$

*where $Y(t)$ is a centered Gaussian process with the covariance*

$$\boldsymbol{E}Y(t)Y(s) = \min(t,s) - ts + \sum_{k \geq 1} (F_k(t,s) + F_k(s,t) - 2ts).$$

**Remark 5.** In the IID case ($\psi(\cdot) \equiv 0$, $b(\cdot, \cdot) \equiv 0$, and $F_k(t,s) = ts$) the conditions of Theorem 2 coincide with the traditional well-known restrictions in Mises, 1947, Filippova, 1959, 1962, Borisov and Sakhanenko, 2000, and others. In this case, the process $Y(t) = W^0(t)$ is a Brownian bridge. Notice also that, in the IID case, the multiple stochastic integral w.r.t the product-noise generated by $W^0(t)$ in Theorem 2 coincides in distribution with the analogous integral w.r.t. a White noise (generated by a standard Wiener process $W(t)$) due to degeneracy of the kernel and the well-known representation $W^0(t) \stackrel{d}{=} W(t) - tW(1)$.

If a kernel $f(t_1, \ldots, t_d)$ vanishes on all diagonal subspaces then sufficient conditions for $f$ to be an element of $S_0$ is as follows: $f \in L_2([0,1]^d, dt_1 \cdots dt_d)$. The kernels of such a kind are used to describe limit behavior of canonical $U$-statistics (for example, see Major, 1994, and the corresponding comment above). Moreover, to describe the Itô – Wiener construction of multiple integrals the kernels vanishing on all diagonal subspaces provide the above-mentioned isometry between closed linear span of multiple integral sums and the corresponding function space.

**Remark 6.** From the beginning, in the IID case, to describe limit behavior of canonical $U$ and $V$-statistics another representation of the limit distribution was used. For instance, in the case $d = 2$, for a canonical $U$-statistics with a symmetric kernel the limit random variable can be represented as the random series

$$\sum_{k=1}^{\infty} \lambda_k (\tau_k^2 - 1), \tag{9}$$

where $\{\tau_k\}$ are i.i.d. standard Gaussian random variables, $\{\lambda_k\}$ are eigenvalues of the integral operator in $L_2[0,1]$ with the kernel $f(s,t)$ (see [25]). Later this result was extended to the general case $d \geq 2$ (see [23]). In the general case the limit random variable can be represented as a polynomial transform of $\{\tau_k\}$ in terms of the Hermite–Chebyshev polynomials. Finally, using analogous representations of multiple Wiener integrals noted in [27], and in [17], one can rewrite the limit random variables as the corresponding multiple stochastic integral w.r.t. to a multiple White noise (see, for example, Dynkin and Mandelbaum, 1983).

In the non-IID case an analog of (9) was obtained in [11], where $\{\lambda_k\}$ are the above-mentioned eigenvalues which are additionally assumed to be summable, and $\{\tau_k\}$ is a Gaussian sequence with the covariance function depending of these eigenvalues as well as of the covariance function of the initial stationary sequence $\{X_i\}$ satisfying $\varphi$-mixing condition.

In the special case when the observations are defined by a nonrandom transform of a Gaussian stationary sequence under another dependency restriction, limit behavior of canonical $U$-statistics were investigated in [8].

*Proof of Theorem 2.* We divide the proof into three steps (for detail see [5]). Let $\{f_k(x_1, \ldots, x_d)\}$ be a sequence of step functions of the form

$$(10) \quad f_k(x_1, \ldots, x_d) := \sum_{j_1, \ldots, j_d \leq k} f_{j_1, \ldots, j_d} \prod_{i=1}^{d} I(x_i \in A_{j_i})$$

converging in $(S_0, ||\cdot||_0)$ to a function $f(x_1, \ldots, x_d)$ as $k \to \infty$. Hence the functions $\{f_k\}$ converge to $f$ in the corresponding linear normed space $(S, ||\cdot||)$ in Theorem 1. Then

S t e p I. For every fixed $k$, due to the multivariate central limit theorem for finite-dimensional distributions of the empirical measure based on slowly dependent observations (see, for example, [1]), as $n \to \infty$,

$$\int_{[0,1]^d} f_k(x_1, \ldots, x_d) dS_n(x_1) \cdots dS_n(x_d) \overset{d}{\to} \int_{[0,1]^d} f_k(x_1, \ldots, x_d) dY(x_1) \cdots dY(x_d).$$

Actually, this weak convergence is valid under weaker dependency conditions than those in Theorem 2.

S t e p II. As $k \to \infty$, by the above-mentioned construction of multiple stochastic integrals (see Theorem 1 and the comment above) we have the following mean square convergence:

$$\int_{[0,1]^d} f_k(x_1, \ldots, x_d) dY(x_1) \cdots dY(x_d) \overset{L_2}{\to} \int_{[0,1]^d} f(x_1, \ldots, x_d) dY(x_1) \cdots dY(x_d).$$

S t e p III. As $k \to \infty$, we should prove the following mean square convergence *uniformly* on $n$:

$$\int_{[0,1]^d} f_k(x_1, \ldots, x_d) dS_n(x_1) \cdots dS_n(x_d) \overset{L_2}{\to} \int_{[0,1]^d} f(x_1, \ldots, x_d) dS_n(x_1) \cdots dS_n(x_d).$$

It is easy to see that the assertion of Theorem 2 follows from items I – III.

So, we should prove only item III. To establish the last convergence, for *every* step function $g_k$ of the form (10), it suffices to prove the estimate

$$(11) \quad \sup_n \mathbf{E} \left( \int_{[0,1]^d} g_k(x_1, \ldots, x_d) dS_n(x_1) \cdots dS_n(x_d) \right)^2 \leq C ||g_k||_0^2,$$

where the constant $C$ depends on $d$ and $\Psi(d)$ only. Indeed, this estimate implies the analogous upper bound for the step function $f_k - f_l$ (see the corresponding comment in the proof of Theorem 1), where $f_k$ and $f_l$ are arbitrary elements of the above-mentioned sequence $\{f_m\}$ converging to the kernel $f$ in the norm $||\cdot||_0$. Therefore, as $l \to \infty$, the corresponding $L_2$-limit of the integral sums $\int f_l(x_1,\ldots,x_d)dS_n(x_1)\cdots dS_n(x_d)$ exists and coincides with the Lebesgue integral in (6) since this limit does not depend on a sequence of step functions approximating to $f$. Thus, as $k \to \infty$,

$$\sup_n \mathbf{E}\left(\int f_k(x_1,\ldots,x_d)dS_n(x_1)\cdots dS_n(x_d) - \int f(x_1,\ldots,x_d)dS_n(x_1)\cdots dS_n(x_d)\right)^2 \to 0$$

for every sequence of step functions $\{f_k\}$ converging to the kernel $f$ in $(S_0, ||\cdot||_0)$ if only relation (11) is valid.

To prove (11) we need the following auxiliary statements.

**Lemma 1.** *(see [21]) Let $\xi$ and $\eta$ be random variables measurable w.r.t. the $\mathcal{F}_1^k$ and $\mathcal{F}_{k+m}^\infty$ $(m \geq 1)$ respectively. If $\mathbf{E}|\xi| < \infty$ and $\mathbf{E}|\eta| < \infty$ then*

(12) $$|\mathbf{E}\xi\eta - \mathbf{E}\xi\mathbf{E}\eta| \leq \psi(m)\mathbf{E}|\xi|\mathbf{E}|\eta|.$$

Introduce the following notation: $\widetilde{I}_k(A) := I(X_k \in A) - P(A)$, $A \in \mathcal{A}$, where the marginal $P$ is the Lebesgue measure on $[0,1]$. The following assertion is easy proved by induction (for detail see [5]).

**Lemma 2.** *For any natural numbers $q$ and $l_1,\ldots,l_q$, and for any pairwise disjoint measurable subsets $A_1,\ldots,A_q$ the following moment inequality is valid:*

(13) $$\mathbf{E}|\widetilde{I}_k^{l_1}(A_1)\cdots\widetilde{I}_k^{l_q}(A_q)| \leq (q+1)P(A_1)\cdots P(A_q).$$

From Lemmas 1 and 2 we deduce the following upper bound.

**Lemma 3.** *Let $q_1 < \cdots < q_s$ be arbitrary natural numbers and $q_0 = 0$. Consider $s$ collections of sets $\{A_1,\ldots A_{q_1}\},\ldots,\{A_{q_{s-1}+1},\ldots A_{q_s}\}$, where the measurable subsets $A_i$ inside of every collection are pairwise disjoint. Put*

$$\nu_{k_i} := \widetilde{I}_{k_i}^{l_{q_{i-1}+1}}(A_{q_{i-1}+1})\cdots\widetilde{I}_{k_i}^{l_{q_i}}(A_{q_i}).$$

*Then, for any natural numbers $k_1 < \cdots < k_s$ and $l_1,\ldots,l_{q_1},\ldots,l_{q_s}$, the following estimate is valid:*

(14) $$\mathbf{E}|\nu_{k_1}\cdots\nu_{k_s}| \leq C(\psi(1),s,q_s)P(A_1)\cdots P(A_{q_s}),$$

*where the constant $C(\cdot)$ depends only on the arguments indicated.*

*Proof.* By (13) we have $\mathbf{E}|\nu_{k_i}| \leq (q_i - q_{i-1} + 1)P(A_{q_{i-1}+1})\cdots\mathbf{P}(A_{q_i})$. It is clear that the random variables $\{\nu_{k_i}\}$ satisfy $\psi$-mixing condition. Hence, by (12) and (13) we obtain

$$\mathbf{E}|\nu_{k_1}\cdots\nu_{k_s}| \leq \prod_{j=1}^{s-1}(1+\psi(k_{j+1}-k_j))\mathbf{E}|\nu_{k_1}|\cdots\mathbf{E}|\nu_{k_s}|$$

$$\leq \prod_{j=1}^{s-1}(1+\psi(k_{j+1}-k_j))\prod_{j=1}^{s-1}(q_j - q_{j-1} + 1)P(A_1)\cdots P(A_{q_s})$$

$$\leq C(\psi(1),s,q_s)P(A_1)\cdots P(A_{q_s}).$$

The Lemma is proved. □

The main assertion to prove (11) is as follows:

**Lemma 4.** *Let $l_1, \ldots, l_q$ be natural numbers such that $l_1 + \cdots + l_q = 2d$, and let $A_1, \ldots, A_q$ be a collection of pairwise disjoint measurable sets. Under condition (8) the following upper bound is valid:*

$$(15) \qquad |\mathbf{E} S_n^{l_1}(A_1) \cdots S_n^{l_q}(A_q)| \leq C P(A_1) \cdots P(A_q).$$

*where the constant $C$ depends only on $d$ and $\Psi(d)$.*

**Remark 7.** The statement of Lemma 4 is also contained in Sen, 1972. However, the corresponding constant in this paper contains the factor $(1 + \psi(0))^2$. It is easy to see that if the marginal distribution has a continuous component (for example, if it is the Lebesgue measure on $[0,1]$) or infinitely many atoms then $\psi(0) = \infty$. To verify this property we may put in (4) $A = B$, where $A$ is an event from the $\sigma$-field $\mathcal{F}_1^1$. So, under the above-mentioned restrictions on the marginal $P$

$$\psi(0) \geq \sup_{A \in \mathcal{F}_1^1, \ P(A) > 0} (1/P(A) - 1) = \infty.$$

In the sequel we will denote by the symbols $C$ or $C_i$ various positive constants depending only on $d$ and $\Psi(d)$.

*Proof of Lemma 4.* We first write the simple estimate

$$|\mathbf{E} S_n^{l_1}(A_1) \cdots S_n^{l_q}(A_q)|$$
$$\leq n^{-d} \sum_{k_1, \ldots, k_{2d} \leq n} |\mathbf{E} \widetilde{I}_{k_1}(A_1) \cdots \widetilde{I}_{k_{l_1}}(A_1) \cdots \widetilde{I}_{k_{2d-l_q+1}}(A_q) \cdots \widetilde{I}_{k_{2d}}(A_q)|.$$

The initial sum on the right-hand side of this inequality can be estimated by a finite sum of the following diagonal subsums

$$(16) \qquad \sum_{k_1 < \cdots < k_r \leq n} |\mathbf{E} \nu_{k_1} \cdots \nu_{k_r}|,$$

where $\nu_{k_i} := \widetilde{I}_{k_i}^{s_1(i)}(A_1) \cdots \widetilde{I}_{k_i}^{s_q(i)}(A_q)$, and the integers $s_j(i)$ are defined by the corresponding diagonal subspace of subscripts in the initial multiple sum and satisfy the conditions $0 \leq s_j(i) \leq l_j$ for all $i \leq r$ and $j \leq q$, and $\sum_{i \leq r} \sum_{j \leq q} s_j(i) = 2d$.

Let $r \leq d$. Estimating by (14) each summand in (16) and taking the normalized factor $n^{-d}$ and the number of summands in (16) into account we obtain the upper bound we need.

Let now $r > d$. We call the random variable $\nu_{k_i}$ *short product* if $\sum_{j \leq q} s_j(i) = 1$, i. e., $\nu_{k_i} = \widetilde{I}_{k_i}(A_{q_i})$ for some $q_i \leq q$. Notice that if $\nu_{k_i}$ is a short product then $\mathbf{E} \nu_{k_i} = 0$.

We now consider the auxiliary multiple sum consisting of the random variables $\nu_{k_i}$ defined in (16) for the fixed diagonal subspace of subscripts:

$$(17) \qquad \sum_{k_{v_1} < \cdots < k_{v_2} \leq n} |\mathbf{E} \nu_{k_{v_1}} \cdots \nu_{k_{v_2}}|,$$

where $1 \leq v_1 < v_2 \leq r$ and the value $v := v_2 - v_1 + 1$ is the dimension of the corresponding multiple sum. Introduce the following notation: $e_j(i) := \min\{1, s_j(i)\}$.

We first prove the following assertion: If, in the summands in (17), there are at least $m$ shorts products, where $0 \leq m \leq v$, then the following upper bound is valid:

$$(18) \quad \sum_{k_{v_1} < \cdots < k_{v_2} \leq n} |\mathbf{E}\nu_{k_{v_1}} \cdots \nu_{k_{v_2}}| \leq Cn^{v-m/2} \prod_{j \leq q} P(A_j)^{\alpha_j(v_1, v_2)},$$

where $\alpha_j(v_1, v_2) := \sum_{i=v_1}^{v_2} e_j(i)$. Notice that the set function $\alpha_j(a, b)$ is additive on intervals $[a, b]$. We prove this assertion by induction on $m$ for all $v_1$ and $v_2$ such that $v \geq m$ and $v \leq r$. Let $m = 1$, i. e., the expectations in (17) contain at least one short product. Denote it by $\nu_{k_l}$, where $k_{v_1} \leq k_l \leq k_{v_2}$. First we note that, in terms of the notation above, we can rewrite the statement of Lemma 3 for the absolute moment of each random product in (17) in such a way:

$$\mathbf{E}|\nu_{k_{v_1}} \cdots \nu_{k_{v_2}}| \leq C \prod_{j \leq q} P(A_j)^{\alpha_j(v_1, v_2)}.$$

Taking this estimate into account we evaluate by (12) every summand in (17) setting $\xi := \nu_{k_{v_1}} \cdots \nu_{k_l}$ and $\eta := \nu_{k_{l+1}} \cdots \nu_{k_{v_2}}$:

$$\sum_{k_{v_1} < \cdots < k_{v_2} \leq n} |\mathbf{E}\nu_{k_{v_1}} \cdots \nu_{k_{v_2}}|$$

$$\leq \sum_{k_{v_1} < \cdots < k_{v_2} \leq n} \psi(k_{l+1} - k_l) \mathbf{E}|\nu_{k_{v_1}} \cdots \nu_{k_l}| \mathbf{E}|\nu_{k_{l+1}} \cdots \nu_{k_{v_2}}|$$

$$+ \sum_{k_{v_1} < \cdots < k_l \leq n} |\mathbf{E}\nu_{k_{v_1}} \cdots \nu_{k_l}| \sum_{k_{l+1} < \cdots < k_{v_2} \leq n} \mathbf{E}|\nu_{k_{l+1}} \cdots \nu_{k_{v_2}}|$$

$$\leq C_1 n^{v-1} \prod_{j \leq q} P(A_j)^{\alpha_j(v_1, v_2)} \sum_{i \geq 1} \psi(i)$$

$$+ C_2 n^{v_2 - l} \prod_{j \leq q} P(A_j)^{\alpha_j(l+1, v_2)}$$

$$\times \sum_{k_{v_1} < \cdots < k_l \leq n} \psi(k_l - k_{l-1}) \mathbf{E}|\nu_{k_1} \cdots \nu_{k_{l-1}}| \mathbf{E}|\nu_{k_l}|$$

$$\leq (C_3 n^{v-1} + C_4 n^{v_2 - l} n^{l - v_1}) \prod_{j \leq q} P(A_j)^{\alpha_j(v_1, v_2)}$$

$$\leq C_5 n^{v - 1/2} \prod_{j \leq q} P(A_j)^{\alpha_j(v_1, v_2)}$$

which required. In this chain of relations the second inequality is valid due to (12) and the equality $\mathbf{E}\nu_{k_l} = 0$ as well.

We now assume that the upper bounds

$$\sum_{k_{v_1} < \cdots < k_{v_2} \leq n} |\mathbf{E}\nu_{k_{v_1}} \cdots \nu_{k_{v_2}}| \leq Cn^{v - z/2} \prod_{j \leq q} P(A_j)^{\alpha_j(v_1, v_2)}$$

are true for all integers $z < m$, where $z$ is the minimal possible number of short products in the expectations under consideration, and for all possible dimensions $v : z \leq v \leq r$ of multiple sums of the form (17), and, moreover, the moments in (17) contain at least $m$ shorts products. Denote these products by $\nu_{k_{l_1}}, \ldots, \nu_{k_{l_m}}$.

Consider the following $m-1$ pairs of neighboring products: $\nu_{k_{l_s}}\nu_{k_{l_s+1}}$, $s = 1\ldots m-1$. Denote by $t_1,\ldots,t_{m-1}$ differences between the subscripts in these pairs. We have

$$\sum_{k_{v_1}<\cdots<k_{v_2}\leq n} |\mathbf{E}\nu_{k_{v_1}}\cdots\nu_{k_{v_2}}| \leq R_1 + \cdots + R_{m-1},$$

where the subsum $R_s$ is taken over the set of subscripts

$$I_s := \{(k_{v_1},\ldots,k_{v_2}) : k_{v_1} < \cdots < k_{v_2} \leq n,\ t_s = \max t_i\}.$$

We now estimate by (12) each summand in $R_s$ setting $\xi := \nu_{k_{v_1}}\cdots\nu_{k_{l_s}}$ and $\eta := \nu_{k_{l_s+1}}\cdots\nu_{k_{v_2}}$:

$$(19) \quad \begin{aligned} R_s &\leq \sum_{I_s} \psi(k_{l_s+1}-k_{l_s})\mathbf{E}|\nu_{k_{v_1}}\cdots\nu_{k_{l_s}}|\mathbf{E}|\nu_{k_{l_s+1}}\cdots\nu_{k_{v_2}}| \\ &+ \sum_{k_{v_1}<\cdots<k_{l_s}\leq n} |\mathbf{E}\nu_{k_{v_1}}\cdots\nu_{k_{l_s}}| \sum_{k_{l_s+1}<\cdots<k_{v_2}\leq n} \mathbf{E}|\nu_{k_{l_s+1}}\cdots\nu_{k_{v_2}}|. \end{aligned}$$

Consider the first sum on the right-hand side of (19):

$$\sum_{I_s} \psi(k_{l_s+1}-k_{l_s})\mathbf{E}|\nu_{k_{v_l}}\cdots\nu_{k_{l_s}}|\mathbf{E}|\nu_{k_{l_s+1}}\cdots\nu_{k_{v_2}}|$$
$$\leq C\prod_{j\leq q} P(A_j)^{\alpha_j(v_1,v_2)} \sum_{I_s} \psi(t_s)$$
$$\leq C\prod_{j\leq q} P(A_j)^{\alpha_j(v_1,v_2)} n^{v-(m-1)} \sum_{t_i:t_i\leq t_s} \psi(t_s)$$
$$\leq C\prod_{j\leq q} P(A_j)^{\alpha_j(v_1,v_2)} n^{v-m/2}\Psi\left(\frac{m}{2}\right).$$

Notice that the last inequality is valid for $m \geq 2$. Consider now the product of the sums on the right-hand side of (19). Let the summands of the first sum contain $m_1$ short products indicated above, and, in the summands of the second sum, there are $m - m_1$ short products indicated above. By the construction we have $1 \leq m_1 \leq m-1$. Hence, for these sums we can use the induction proposition. Finally, we have

$$\sum_{k_{v_1}<\cdots<k_{l_s}\leq n} |\mathbf{E}\nu_{k_{v_1}}\cdots\nu_{k_{l_s}}| \sum_{k_{l_s+1}<\cdots<k_{v_2}\leq n} |\mathbf{E}\nu_{k_{l_s+1}}\cdots\nu_{k_{v_2}}|$$
$$\leq Cn^{l_s-m_1/2}n^{v-l_s-(m-m_1)/2}\prod_{j\leq q} P(A_j)^{\alpha_j(v_1,v_2)}$$
$$= Cn^{v-m/2}\prod_{j\leq q} P(A_j)^{\alpha_j(v_1,v_2)}$$

which required. Thus, for $R_s$ we obtained the upper bound we need. It means that the analogous estimate is valid for whole sum (17). The induction is over.

To prove the assertion of Lemma 4 we should note that, first, by the definition, $\alpha_j(1,r) \geq 1$ for all $j \leq q$, and, second, in the case $r > d$, the summands in (16) contain at least $2(r-d)$ short products. So, we should put in (18) $v_1 := 1$, $v_2 := r$, $m := 2(r-d)$ and $v := r$. It means that, for the sum in (16), the following upper bound is valid:

$$\sum_{k_1<\cdots<k_r\leq n} |\mathbf{E}\nu_{k_1}\cdots\nu_{k_r}| \leq Cn^d \prod_{j\leq q} P(A_j)^{\alpha_j(1,r)} \leq Cn^d P(A_1)\cdots P(A_q).$$

The Lemma is proved. □

In conclusion we deduce from Lemma 4 the upper bound (11). For every step function $g_k$ of the form (10) we have

$$
\begin{aligned}
\sup_n \mathbf{E} & \left( \int_{[0,1]^d} g_k(x_1, \ldots, x_d) dS_n(x_1) \cdots dS_n(x_d) \right)^2 \\
& = \sup_n \mathbf{E} \int_{[0,1]^{2d}} g_k(x_1, \ldots x_d) g_k(x_{d+1}, \ldots x_{2d}) dS_n(x_1) \cdots dS_n(x_{2d}) \\
& \leq \sup_n \sum_{i_1, \ldots, i_{2d}} |g_{i_1, \ldots, i_d}||g_{i_{d+1}, \ldots, i_{2d}}| \mathbf{E}|S_n(A_{i_1}) \cdots S_n(A_{i_{2d}})| \\
& \leq C \sum_{i \geq 0} \int_{D_i^{(2)}} |g_k(x_1, \ldots x_d)||g_k(x_{d+1}, \ldots x_{2d})| P(dx_{i_{q_1}}) \cdots P(dx_{i_{q_r}}),
\end{aligned}
\tag{20}
$$

where $D_i^{(2)}$ is a diagonal (or the main) subspace in $[0,1]^{2d}$ of dimension $r$ which is defined by the integers $q_1, \ldots, q_r$ by analogy with $D_i$ in $[0,1]^d$. Further, by the Cauchy – Bunyakovskii inequality we estimate every multiple integral on the right-hand side of (20):

$$
\begin{aligned}
\int_{D_i^{(2)}} & |g_k(x_1, \ldots x_d)||g_k(x_{d+1}, \ldots x_{2d})| P(dx_{i_{q_1}}) \cdots P(dx_{i_{q_r}}) \\
& \leq \left( \int_{D_{j_1}} g_k^2(x_1, \ldots x_d) P(dx_{i_{v_1}}) \cdots P(dx_{i_{v_l}}) \right. \\
& \quad \times \left. \int_{D_{j_2}} g_k^2(x_1, \ldots x_d) P(dx_{i_{s_1}}) \cdots P(dx_{i_{s_m}}) \right)^{1/2} \\
& \leq \|g_k\|_0^2,
\end{aligned}
\tag{21}
$$

where $D_{j_1}$ and $D_{j_2}$ are the corresponding diagonal (or the main) subspaces in $[0,1]^d$ of respective dimensions $l$ and $m$, where $l + m \geq r$, defined by the integers $\{v_j\}$ and $\{s_j\}$. Thus, the upper bound in (11) is proved. □